\documentclass[12pt]{amsart}
\usepackage{amssymb}

\theoremstyle{plain}
\newtheorem*{theorem}{Theorem}

\theoremstyle{definition}
\newtheorem{example}{Example}
\newtheorem*{remark}{Remark}

\newcommand{\mat}{\begin{bmatrix}a&b\\c&d\end{bmatrix}}

\newcommand{\gie}{\mathbb{G}}

\renewcommand{\geq}{\geqslant}

\begin{document}
\title[Eigenvectors of semigroups on von Neumann algebras]
{The eigenvectors of semigroups \\
of positive maps on von Neumann algebras}
\author[A. \L uczak]{Andrzej \L uczak}
\address{Faculty of Mathematics and Computer Science\\
        \L\'od\'z University \\
        ul. S. Banacha 22 \\
        90-238 \L\'od\'z, Poland}
\email{anluczak@math.uni.lodz.pl}
\thanks{Work supported by KBN grant 2PO3A 03024}
\date{}
\begin{abstract}
The eigenvectors of an ergodic semigroup of linear normal positive
unital maps on a von Neumann algebra are described. Moreover, it
is shown by means of examples, that mere positivity of the maps in
question is not sufficient for Frobenius theory as in \cite{A-H}
to hold.
\end{abstract}
\maketitle
\section{Introduction}\label{S0}
Frobenius theory for completely positive maps on von Neumann
algebras was developed in \cite{A-H} (further contributions to
this subject can be found in \cite{Gr} and \cite{Wa}). This theory
states, in particular, that for an ergodic semigroup of completely
positive (or, in fact, even Schwarz) maps on a von Neumann algebra
its point spectrum forms a group, and the corresponding
eigenspaces are one-dimensional and spanned by a unitary operator.
The aim of this paper is to investigate the structure of the point spectrum of a semigroup of positive maps on a von Neumann algebra. Namely, we prove that the eigenvectors are either multiples of a partial isometry or linear combinations of two partial isometries or multiples of a unitary operator from this algebra. As for Frobenius theory, we show by means of examples, that there is a whole class of semigroups of positive maps on a von Neumann algebra such that their point spectra are not groups and the corresponding eigenspaces have dimensions greater than one.
\section{Preliminaries and notation}\label{S1}
Let $M$ be a von Neumann algebra, and let $(\varPhi_g\colon
g\in\gie)$ be a semigroup of linear normal positive unital maps on
$M$. We shall be concerned with two cases: $\gie=\mathbb{N}_0$ --- all nonnegative integers, and \linebreak $\gie=\mathbb{R}_+$ --- all
nonnegative reals (notice that in the first case the semigroup has
the form $(\varPhi^n: n=0,1,\dots)$, where $\varPhi$ is a
linear normal positive unital map on $M$).

A complex number $\lambda$ of modulus one is called an
\emph{eigenvalue} of the semigroup, if there is a nonzero $x\in M$
such that for each $g\in\gie$
\begin{equation}\label{e0}
 \varPhi_g(x) = \lambda^g x.
\end{equation}
The collection of all $x's$ such that \eqref{e0} holds is called
the \emph{eigenspace} corresponding to the eigenvalue $\lambda$,
and denoted by $M_{\lambda}$. In particular, $M_1$ is the
fixed-point space of the semigroup, and the semigroup is called
\emph{ergodic} if $M_1$ consists of multiples of the identity. The
set of all eigenvalues of the semigroup is called its \emph{point
spectrum}, and denoted by $\sigma((\varPhi_g))$. Let $\omega$ be a
normal faithful state on $M$ such that for each \linebreak
$g\in\gie,\ \omega\circ\varPhi_g = \omega$. The part of Frobenius
theory developed in \cite{A-H} which is of interest to us, states
that in this case, if we assume that \linebreak $(\varPhi_g: g\in\gie)$
is ergodic and the maps $\varPhi_g$ are two-positive, the
point spectrum is a group, and the corresponding eigenspaces are
one-dimensional and spanned by a unitary operator. A natural question
is if the same is true under the assumption of mere positivity of the
maps $\varPhi_g$. We shall show that this is not the case
for the group structure nor for the dimension of the eigenspaces.

Let $N$ be the $\sigma$-weak closure of the linear span of
$\bigcup_{\lambda}M_{\lambda}$, where the sum is taken over all
eigenvalues of $(\varPhi_g)$. Then according to \linebreak \cite[Theorem
1]{Lu} $N$ is a $JW^*$-algebra, by which is meant that $N$ is
\linebreak a $\sigma$-weakly closed linear space, closed with respect
to the Jordan product
\[
 x\circ y = \frac{1}{2}(xy + yx);
\]
moreover, $\varPhi_g|N$ are Jordan $^*$-automorphisms.
Consequently, if $x\in M_{\lambda}$, then
\[
 \varPhi_g(x^*) = \varPhi_g(x)^* = \bar{\lambda}^g x^*,
\]
meaning that $x^*\in M_{\bar{\lambda}}$, and for $x\in
M_{\lambda_1},\ y\in M_{\lambda_2}$ we have
\[
 \varPhi_g(x\circ y) = \varPhi_g(x)\circ\varPhi_g(y) = \lambda_1^g
 x\circ\lambda_2^g y = (\lambda_1\lambda_2)^g x\circ y,
\]
meaning that $x\circ y\in M_{\lambda_1\lambda_2}$. In particular,
for an eigenvector $x$ we have $x\circ x^*\in M_1$.

\section{Main theorem and examples}\label{S2}
In the following theorem we describe the eigenvectors of
$(\varPhi_g)$.
\begin{theorem}
Assume that $(\varPhi_g)$ is ergodic, and let $x\in M_{\lambda}$
be an eigenvector of $(\varPhi_g)$. Then one
of the following possibilities holds:
\begin{enumerate}
\item[(i)] $x=\alpha v$, where $\alpha\in\mathbb{C}$, and $v$ is a
partial isometry in $M_{\lambda}$ such that
\[
 v^*v=e, \qquad vv^*=e^{\bot}
\]
for some nonzero projection $e$, $e\neq\boldsymbol{1}$;
\item[(ii)] $x=\alpha_1v_1+\alpha_2v_2$, where
$\alpha_1,\alpha_2$ are nonzero complex numbers, $\alpha_1\neq\alpha_2$, and
$v_1,v_2$ are partial isometries in $M_{\lambda}$ such that for
some nonzero projection $e,\ e\neq\boldsymbol{1}$,
\[
 v^*_1v_1=e, \quad v_1v^*_1=e^{\bot}, \quad v^*_2v_2=e^{\bot}, \quad
 v_2v^*_2=e;
\]
\item[(iii)] $x=\alpha u$, where $\alpha\in\mathbb{C}$, and $u$ is
a unitary operator in $M_{\lambda}$.
\end{enumerate}
\end{theorem}

Let us begin with a simple example.
\begin{example}\label{Ex1}
Put $M=\mathbb{B}(\mathbb{C}^2)$, $\omega=\frac{1}{2}\operatorname{tr}$, and let $\lambda_0\in\mathbb{C}$ be such
that $|\lambda_0|=1,\ \lambda_0\ne1$. Define
$\varPhi\colon M\to M$ as
\[
 \varPhi\left(\mat\right) =
 \begin{bmatrix}\frac{a+d}{2}&\lambda_0b\\
 \bar{\lambda}_0c&\frac{a+d}{2}\end{bmatrix}.
\]
$\varPhi$ is linear normal unital positive, and $\omega\circ\varPhi=\omega$. Moreover, $(\varPhi^n)$ is
ergodic.

(i) Let $\lambda_0$ be such that $\lambda_0\ne-1,\
\lambda_0^3\ne1$. Then it turns out that $\sigma((\varPhi^n)) =
\{1,\lambda_0,\bar{\lambda}_0\}$, which is not a group if $\lambda_0\ne-1,\ \lambda_0^3\ne1$.

(ii) Now let $\lambda_0 = -1$. Then we obtain
$\sigma((\varPhi^n)) = \{1,-1\}$, and
\[
 M_{-1} = \left\{\begin{bmatrix} 0 & b \\ c & 0
 \end{bmatrix}: b,c\in\mathbb{C}\right\},
\]
so the eigenspace is not one-dimensional.\hfill\qedsymbol
\end{example}

The same conclusion may be obtained for the semigroup $(\varPhi_t: t\geq 0)$ defined as
\[
 \varPhi_t\left(\mat\right) = \begin{bmatrix} \frac{a+d}{2} &
 \lambda_0^t b \\ \bar{\lambda}_0^t c & \frac{a+d}{2}
 \end{bmatrix},
\]
thus giving a corresponding example in the continuous case.
\begin{remark}
The triple $(\mathbb{B}(\mathbb{C}^2),(\varPhi_t),\omega)$ from
the above example constitutes what in \cite{Gr} is called an
irreducible $W^*$-dynamical system. In \linebreak\cite[Theorem
3.8]{Gr} it is proved that under the assumption that the
$\varPhi_t$'s are Schwarz maps every such system on the full
algebra has trivial point spectrum (i.e. consisting only of $1$).
As we see this is not the case if we assume only positivity of the
maps $\varPhi_t$'s.
\end{remark}
Now we construct a more involved
example (in fact, a class of examples) in which we shall see that
all the possibilities for the eigenvectors given in the Theorem may
occur.
\begin{example}\label{Ex2}
Let $M$ be abelian, let $\omega$ be a normal faithful state on
$M$, and let $\varPsi$ be a positive normal unital map on $M$ such
that $\omega\circ\varPsi = \omega,\ (\varPsi^n)$ is ergodic, and
$\sigma((\varPsi^n)) = \{-1,1\}$. The abelianess of $M$ implies
that $\varPsi$ is completely positive (cf. \cite[Chapter IV.3]{Ta}),
thus according to \cite{A-H} the eigenspace corresponding to $-1$
is one-dimensional and spanned by a unitary operator $u$, i.e.
\[
 \varPsi(x) = -x
\]
if and only if $x$ is a multiple of $u$.

Put $\widetilde{M} = \textbf{Mat}_2(M)$ --- the algebra of all $2\times 2$-matrices with elements from $M$,
\[
 \widetilde\omega\left(\mat\right) = \frac{1}{2}[\omega(a) +
 \omega(d)],
\]
and let $\lambda_0\in\mathbb{C}$ be such that $|\lambda_0|=1,\
\lambda_0\notin\{-1,1\}$. Define
$\varPhi\colon\widetilde{M}\to\widetilde{M}$ as
\[
 \varPhi\left(\mat\right) =
 \begin{bmatrix}\varPsi\left(\frac{a+d}{2}\right) & \lambda_0\varPsi(b) \\
 \bar{\lambda}_0\varPsi(c) & \varPsi\left(\frac{a+d}{2}\right)\end{bmatrix}.
\]
$\varPhi$ is a linear normal unital positive map on $\widetilde{M}$ such
that $\widetilde\omega\circ\varPhi = \widetilde\omega$. Moreover, $(\varPhi^n)$ is ergodic and
\[
 \widetilde{M}_1=\left\{\begin{bmatrix}
    \theta\boldsymbol{1} & 0 \\
    0 & \theta\boldsymbol{1}
    \end{bmatrix}:\theta\in\mathbb{C}\right\}.
\]

Let $\lambda\ne1$ be an eigenvalue of $\varPhi$.

(i) Take $\lambda_0$ such that $\lambda_0\neq i,\
\lambda_0\neq-i,\ \lambda^3_0\neq1,\ \lambda^3_0\neq-1$.
We have the following possibilities:

(i.1) $\lambda=-1$. Then
\[
 \widetilde{M}_{-1}=\left\{\alpha\begin{bmatrix}u & 0\\ 0 & u
 \end{bmatrix}: \alpha\in\mathbb{C}\right\},
\]
so the eigenvector corresponding to the eigenvalue $-1$ is as in
part (iii) of the Theorem.

(i.2) $\lambda\neq-1$. In this case one of the four
situations must occur:

(i.2.1) $\lambda=\lambda_0$. Then
\[
 \widetilde{M}_{\lambda_0}=\left\{\alpha\begin{bmatrix}0 &
 \boldsymbol{1}\\ 0 & 0
 \end{bmatrix}:\alpha\in\mathbb{C}\right\},
\]
so the eigenvector corresponding to the eigenvalue $\lambda_0$ is
as in part (i) of the Theorem.

(i.2.2) $\lambda=-\lambda_0$. Then
\[
 \widetilde{M}_{-\lambda_0}=\left\{\alpha\begin{bmatrix}0 & u\\
 0 & 0\end{bmatrix}:\alpha\in\mathbb{C}\right\}.
\]

(i.2.3) $\lambda=\bar{\lambda}_0$. Then
\[
 \widetilde{M}_{\bar{\lambda}_0}=\left\{\alpha\begin{bmatrix}0 &
 0\\ \boldsymbol{1} & 0
 \end{bmatrix}:\alpha\in\mathbb{C}\right\}.
\]

(i.2.4) $\lambda=-\bar{\lambda}_0$. Then
\[
 \widetilde{M}_{-\bar{\lambda}_0}=\left\{\alpha\begin{bmatrix}0 &
 0\\u & 0 \end{bmatrix}:\alpha\in\mathbb{C}\right\}.
\]

Moreover, we have
\[
 \sigma((\varPhi^n))=\{1,-1,\lambda_0,\bar{\lambda}_0,-\lambda_0,
 -\bar{\lambda}_0\},
\]
which is not a group under our assumptions on $\lambda_0$.

(ii) Now take $\lambda_0=i$. As in part (i) we have the possibilities:

(ii.1) $\lambda=-1$. Then
\[
 \widetilde{M}_{-1}=\left\{\begin{bmatrix}
    \theta\boldsymbol{1} & 0 \\
    0 & \theta\boldsymbol{1}
    \end{bmatrix}:\theta\in\mathbb{C}\right\}.
\]

(ii.2) $\lambda\neq-1$. In this case we may only
have either $\lambda=i$ or $\lambda=-i$. In the first case
\[
 \widetilde{M}_i=\left\{\alpha_1\begin{bmatrix}0 &
 \boldsymbol{1}\\ 0 & 0\end{bmatrix}+ \alpha_2\begin{bmatrix}0 &
 0\\ u &
 0\end{bmatrix}:\alpha_1,\alpha_2\in\mathbb{C}\right\},
\]
so the situation is as in part (ii) of the Theorem with
\[
 v_1=\begin{bmatrix}0 & \boldsymbol{1}\\ 0 & 0\end{bmatrix},\
 v_2=\begin{bmatrix}0 & 0\\ u & 0\end{bmatrix},\
 e=\begin{bmatrix}0 & 0\\ 0 & \boldsymbol{1}\end{bmatrix},\
 e^{\bot}=\begin{bmatrix}\boldsymbol{1} & 0\\ 0 & 0\end{bmatrix}.
\]
In the second case
\[
 \widetilde{M}_{-i}=\left\{\alpha_1\begin{bmatrix}0 & u\\ 0 &
 0\end{bmatrix}+ \alpha_2\begin{bmatrix}0 & 0\\ \boldsymbol{1} &
 0\end{bmatrix}:\alpha_1,\alpha_2\in\mathbb{C}\right\},
\]
so again part (ii) of the Theorem occurs with
\[
 v_1=\begin{bmatrix}0 & u\\ 0 & 0\end{bmatrix},\
 v_2=\begin{bmatrix}0 & 0\\ \boldsymbol{1} & 0\end{bmatrix},\
 e=\begin{bmatrix}\boldsymbol{1} & 0\\ 0 & 0\end{bmatrix},\
 e^{\bot}=\begin{bmatrix}0 & 0\\ 0 & \boldsymbol{1}\end{bmatrix}.
\]
\hfill\qedsymbol
\end{example}
As in Example \ref{Ex1}, we observe that taking $(\varPsi_t:
t\geq 0)$---an ergodic semigroup of positive maps on $M$ with
$\sigma((\varPsi_t)) = \{1,-1\}$, and defining
\[
 \varPhi_t\left(\mat\right) =
 \begin{bmatrix}\varPsi_t\left(\frac{a+d}{2}\right) &
 \lambda_0^t\varPsi_t(b) \\ \bar{\lambda}_0^t\varPsi_t(c) &
 \varPsi_t\left(\frac{a+d}{2}\right)\end{bmatrix},\qquad t\geq 0,
\]
we obtain a continuous counterpart of Example \ref{Ex2}.


\begin{thebibliography}{9}
 \bibitem{A-H}
  S. Albeverio and R. H\o egh-Krohn, \emph{Frobenius theory of
  positive maps of von Neumann algebras}, Comm. Math. Phys.
  \textbf{64}(1978), 83--94.
 \bibitem{Gr}
  U. Groh, \emph{Positive semigroups on $C^*$- and $W^*$-algebras}, in
  ''One-parameter Semigroups of Positive Operators'', Lecture
  Notes in Math. \textbf{1184}(1986), 369--425.
 \bibitem{Lu}
  A. \L uczak, \emph{Eigenvalues and eigenspaces of quantum
  dynamical systems and their tensor products}, J. Math. Anal.
  Appl. \textbf{221}(1998), 13--32.
 \bibitem{Ta}
  M. Takesaki, \emph{Theory of Operator Algebras}, Springer,
  Berlin--Heidelberg--New York, 1979.
 \bibitem{Wa}
  S. Watanabe, \emph{Asymptotic behavior and eigenvalues of
  dynamical semi-groups on operator algebras}, J. Math. Anal.
  Appl. \textbf{86}(1982), 411--424.
\end{thebibliography}
\end{document}